\documentclass[11pt]{article}
\usepackage{amsmath}
\usepackage{amssymb}
\usepackage{amsfonts}
\usepackage{mathrsfs}
\usepackage{shortvrb,psfrag}
\usepackage{enumerate}
\usepackage{color}
\usepackage[dvips]{graphicx}

\let\al=\alpha
\let\be=\beta

\let\f=\frac
\let\p=\psi

\let\Om=\Omega

\def\R{\Bbb R}

\def\no{\noindent}
\def\na{\nabla}
\def\p{\partial}
\def\eqdef{\buildrel\hbox{\footnotesize def}\over =}
\def\endproof{\hphantom{MM}\hfill\llap{$\square$}\goodbreak}

\newcommand{\beq}{\begin{equation}}
\newcommand{\eeq}{\end{equation}}
\newcommand{\ben}{\begin{eqnarray}}
\newcommand{\een}{\end{eqnarray}}
\newcommand{\beno}{\begin{eqnarray*}}
\newcommand{\eeno}{\end{eqnarray*}}

\newtheorem{Theorem}{Theorem}[section]
\newtheorem{Definition}[Theorem]{Definition}
\newtheorem{Proposition}[Theorem]{Proposition}
\newtheorem{Lemma}[Theorem]{Lemma}

\newtheorem{Remark}[Theorem]{Remark}

\begin{document}
\title{On the interior regularity criterion and the number of singular points to the  Navier-Stokes equations}

\author{Wendong Wang \,and\, Zhifei Zhang\\[2mm]
{\small School of  Mathematical Sciences and BICMR, Peking University, Beijing 100871, P.R. China}\\[1mm]
{\small E-mail: wendong@math.pku.edu.cn and
zfzhang@math.pku.edu.cn}}

\date{January 3, 2012}
\maketitle

\begin{abstract}
We establish some interior regularity criterions of suitable weak
solutions for the 3-D Navier-Stokes equations, which allow the vertical part of the velocity to be large
under the local scaling invariant norm. As an application, we improve Ladyzhenskaya-Prodi-Serrin's criterion and
Escauriza-Seregin-\v{S}ver\'{a}k's criterion. We also show that if weak solution $u$ satisfies
\beno
\|u(\cdot,t)\|_{L^p}\le C(-t)^{\f {3-p}{2p}}
\eeno
for some $3<p<\infty$, then the number of singular points is finite.
\end{abstract}

\setcounter{equation}{0}
\section{Introduction}
We consider the three dimensional incompressible Navier-Stokes equations
\begin{equation}\label{eq:NS}
\left\{\begin{array}{l}
\partial_t u-\Delta u+u\cdot \nabla u+\nabla \pi=0,\\
{\rm div } u=0,
\end{array}\right.
\end{equation}
where $u(x,t)=(u_1(x,t),u_2(x,t),u_3(x,t))$ denotes the unknown
velocity of the fluid, and the scalar function $\pi(x,t)$ denotes the
unknown pressure.

In a seminal paper \cite{Leray}, Leray proved the global existence of weak solution with finite energy.
It is well known that weak solution is unique and regular in two spatial dimensions.
In three dimensions, however, the question of regularity and uniqueness of weak solution is an
outstanding open problem in mathematical fluid mechanics.

In a fundamental paper \cite{CKN},
Caffarelli-Kohn-Nirenberg proved that one-dimensional Hausdorff measure
of the possible singular points of suitable weak solution $u$ is zero (see also \cite{Lin, Tian, LS, Vasseur}).
The proof is based on the following $\varepsilon$-regularity criterion: there exists an $\varepsilon>0$ such that if $u$ satisfies
\ben\label{eq:CKN}
\limsup_{r\rightarrow 0}r^{-1}\int_{Q_r(z_0)}|\nabla u(y,s)|^2dyds\le \varepsilon,
\een
then $u$ is regular at $z_0$. The same result remains true if (\ref{eq:CKN}) is replaced by
\ben\label{eq:CKN1}
\limsup_{r\rightarrow 0}r^{-2}\int_{Q_r(z_0)}|u(y,s)|^3dyds\leq \varepsilon.
\een
The quantities on the left hand side of (\ref{eq:CKN}) and (\ref{eq:CKN1}) are scaling invariant.
More general interior regularity criterions were obtained by Gustafson-Kang-Tsai \cite{GKT}
in terms of scaling invariant quantities (see Proposition \ref{prop:small regularity-GKT}).
In the first part of this paper, we will establish some interior regularity criterions,  which allow the vertical part of the velocity to be large
under the local scaling invariant norm.
The proof is based on the blow-up argument and an observation that if the horizontal part of the velocity is small, then
the blow-up limit satisfies $u_h=0$, hence $\p_3u_3=0$ and
\beno
\p_t u_3-\Delta u_3+\p_3\pi=0,\quad \Delta \pi=0.
\eeno
Using new interior regularity criterions, we improve Ladyzhenskaya-Prodi-Serrin regularity criterions,
which state if the weak solution $u$ satisfies
\beno
u\in L^q(0,T;L^p(\R^3))\quad\textrm{ with} \quad \f 2 q+\f 3p\le 1,\, p\ge 3,
\eeno
then it is regular in $(0,T)\times \R^3$, see \cite{Serrin, Giga, Struwe, ESS}. It should be pointed out
that the regularity in the class $L^\infty(0,T;L^3(\R^3))$ is highly nontrivial,
since it does not fall in the framework of small energy regularity.
This case was solved by Escauriza-Seregin-\v{S}ver\'{a}k \cite{ESS} by using blow-up analysis and the
backward uniqueness for the parabolic equation.

In Leary's paper \cite{Leray}, he also proved that if $[-T,0)$ is the maximal existence interval of smooth solution, then for $p>3$, there exits
$c_p>0$ such that
\beno
\|u(\cdot,t)\|_{L^p}\ge c_p(-t)^{\f {3-p}{2p}}.
\eeno
In general, if $u$ satisfies
\ben\label{eq:blow-up}
\|u(\cdot,t)\|_{L^p}\le C(-t)^{\f {3-p}{2p}},
\een
the regularity of the solution at $t=0$ remains unknown except $p=3$.
Recently, for the axisymmetric Navier-Stokes equations, important progress has been made by
Chen-Strain-Yau-Tsai \cite{chen2, chen} and Koch-Nadirashvili-Segegin-\v{S}ver\'{a}k \cite{KNSS},
where they showed that the solution does not develop Type I singularity (i.e, $\|u(\cdot,t)\|_{L^\infty}\le C(-t)^{-\f12}$)
by using De-Giorgi-Nash method and Liouville theorem respectively.
However, the case without the axisymmetric assumption is still open.
The second part of this paper will be devoted to show that
the number of singular points is finite if the solution satisfies (\ref{eq:blow-up}) for $3<p<\infty$.
The proof is based on an improved $\varepsilon$-regularity criterion: if the suitable weak solution $(u,\pi)$
satisfies
\beno
&&\sup_{t\in [-1+t_0,t_0]}\int_{B_1(x_0)}|u(x,t)|^2dx+\int_{-1+t_0}^{t_0}\big(\int_{B_1(x_0)}|u(x,t)|^4dx\big)^{\f12}dt\\
&&\quad+\int_{-1+t_0}^{t_0}\big(\int_{B_1(x_0)}|\pi(x,t)|^2dx\big)^{\f12}dt \leq \varepsilon_6,
\eeno
then $u$ is regular in $Q_{\f1 2}(z_0)$, see Proposition \ref{prop:small regularity-new}.

This paper is organized as follows. In section 2, we introduce some definitions and notations.
In section 3, we establish some new interior regularity criterions of suitable weak solutions.
In section 4, we apply them to improve Ladyzhenskaya-Prodi-Serrin's criterion and Escauriza-Seregin-\v{S}ver\'{a}k's criterion.
Section 4 is devoted to the proof of the number of singular points under the condition (\ref{eq:blow-up}).
In the appendix, we present the estimates of the pressure and some scaling invariant quantities.

\section{Definitions and notations}

Let us first introduce the definition of weak solution.

\begin{Definition} Let $\Om\subset \R^3$ and $T>0$. We say that $u$ is a Leray-Hopf weak solution
of (\ref{eq:NS}) in $\Om_T=\Om\times (-T,0)$ if
\begin{enumerate}

\item $u\in L^{\infty}(-T,0;L^2(\Om))\cap L^2(-T,0;H^1(\Om))$;

\item $u$ satisfies (\ref{eq:NS}) in the sense of distribution;

\item $u$ satisfies the energy inequality: for a.e. $t\in
[-T,0]$,
\beno \int_{\Om}|u(x,t)|^2dx+2\int_{-T}^t\int_{\Om}|\nabla
u|^2 dxds\leq \int_{\Om}|u(x,-T)|^2dx.
\eeno
\end{enumerate}
Furthermore, the pair $(u,\pi)$ is called a suitable weak solution if $\pi\in L^{3/2}(\Om_T)$
and the energy inequality is replaced by the following local energy inequality:
for  any nonnegative $\phi\in C_c^\infty(\R^3\times\R)$
vanishing in a neighborhood of the parabolic boundary of $\Om_T$,
\beno
&&\int_{\Om}|u(x,t)|^2\phi dx+2\int_{-T}^t\int_{\Om}|\nabla u|^2\phi dxds\\
&&\quad\leq
\int_{-T}^t\int_{\Om}|u|^2(\partial_s\phi+\triangle\phi)+u\cdot\nabla\phi(|u|^2+2\pi)dxds,\quad \textrm{for a.e. } t\in [-T,0].
\eeno
\end{Definition}

\begin{Remark}\label{rem:weak solution}
In general, we don't know whether a Leray-Hopf weak solution is a
suitable weak solution. However, if $u$ is a Leray-Hopf weak
solution and $u\in L^4(\Om_T)$, then it is also a suitable weak
solution, which can be verified by using a standard mollification procedure.
\end{Remark}

Let $(u,\pi)$ be a solution of (\ref{eq:NS}) and introduce the following scaling
\ben\label{eq:scaling}
u^{\lambda}(x, t)={\lambda}u(\lambda x,\lambda^2 t),\quad \pi^{\lambda}(x, t)={\lambda}^2\pi(\lambda x,\lambda^2 t),
\een
for any $\lambda> 0,$ then the family $(u^{\lambda}, \pi^{\lambda})$ is also a solution of (\ref{eq:NS}).
We introduce some invariant quantities under the scaling (\ref{eq:scaling}):
\beno
&&A(u,r,z_0)=\sup_{-r^2+t_0\leq t<t_0}r^{-1}\int_{B_r(x_0)}|u(y,t)|^2dy,\\
&&C(u,r,z_0)=r^{-2}\int_{Q_r(z_0)}|u(y,s)|^3dyds,\\
&&E(u,r,z_0)=r^{-1}\int_{Q_r(z_0)}|\nabla u(y,s)|^2dyds,\\
&&D(\pi,r,z_0)=r^{-2}\int_{Q_r(z_0)}|\pi(y,s)|^{\f32}dyds,
\eeno
where $z_0=(x_0,t), Q_r(z_0)=(-r^2+t_0,t_0)\times B_r(x_0)$,
and $B_r(x_0)$ is a ball of radius $r$ centered at $x_0$. We also
denote $Q_r$ by $Q_r(0)$ and $B_r$ by $B_r(0)$. We also denote
\beno
&&G(f,p,q;r,z_0)=r^{1-\frac3p-\frac2q}\|f\|_{L^{p,q}(Q_r(z_0))},\\
&&H(f,p,q;r,z_0)=r^{2-\frac3p-\frac2q}\|f\|_{L^{p,q}(Q_r(z_0))},\\
&&\widetilde{G}(f,p,q;r,z_0)=r^{1-\frac3p-\frac2q}\|f-(f)_{B_r(x_0)}\|_{L^{p,q}(Q_r(z_0))},\\
&&\widetilde{H}(f,p,q;r,z_0)=r^{2-\frac3p-\frac2q}\|f-(f)_{B_r(x_0)}\|_{L^{p,q}(Q_r(z_0))},
\eeno
where the mixed space-time norm $\|\cdot\|_{L^{p,q}(Q_r(z_0))}$ is defined by
\beno
\|f\|_{L^{p,q}(Q_r(z_0))}\eqdef \Big(\int_{t_0-r^2}^{t_0}\Big(\int_{B_r(x_0)}|f(x,t)|^pdx\Big)^{\f
q p}dt\Big)^\f 1q,
\eeno
and $(f)_{B_r(x_0)}$ is the average of $f$ in the ball $B_r(x_0)$. For the simplicity of notations, we denote
$$A(u,r,(0,0))=A(u,r),\quad \tilde{C}(u,r)=C(u-(u)_{B_r},r),\quad G(f,p,q;r,(0,0))=G(f,p,q;r)$$
and so on. These scaling invariant quantities will play an important role in the interior regularity theory.

Now we recall the definitions of Lorentz space and BMO space \cite{Graf}.

\begin{Definition} Let $\Omega\subset\R^n$ and $1\leq p,\ell \leq \infty$.
We say that a measurable function $f\in L^{p,\ell}(\Omega)$ if $\|f\|_{L^{p,\ell}(\Omega)}<+\infty$, where
\beno
\|f\|_{L^{p,\ell}(\Omega)}\eqdef
\left\{\begin{array}{l}
\Big(\int_0^{\infty}\sigma^{\ell-1}|\{x\in \Omega;|f|>\sigma\}|^{\frac{\ell}{p}}d\sigma\Big)^{\f 1{\ell}}\quad
\textrm{for  } \ell<+\infty,\\
\displaystyle\sup_{\sigma>0}\sigma|\{x\in \Omega;|f|>\sigma\}|^{\frac{1}{p}}\quad
\textrm{for  } \ell=+\infty.
\end{array}\right.
\eeno
Moreover, $f(x,t)\in L^{q,s}(-T,0; L^{p,\ell}(\Omega))$ if $\|f(\cdot,t)\|_{L^{p,\ell}(\Omega)}\in L^{q,s}(-T,0)$.
\end{Definition}
The following facts will be used frequently: for any $R>0$,
\ben\label{eq:lorentz-inc}
&&\|f\|_{L^{p,\ell_1}}\le \|f\|_{L^{p,\ell_2}},\quad\textrm{ if }\quad \ell_1\ge \ell_2;\\
&&\|f\|_{L^{p_1}(\Om)}^{p_1}\le C\big(R^{p_1}|\Om|+R^{p_1-p}\|f\|_{L^{p,\infty}(\Om)}^p\big),\quad\textrm{ if }p>p_1.\label{eq:lorentz}
\een

Recall that a local integrable function $f\in BMO(\R^n)$ if it satisfies
\beno
\sup_{R>0,x_0\in\R^n}\frac{1}{|B_R(x_0)|}\int_{B_R(x_0)}|f(x)-f_{B_R(x_0)}|dx<\infty.
\eeno
Moreover, $f(x)\in VMO(\R^n)$ if $f(x)\in BMO(\R^n)$ and for any $x_0\in \R^n$,
$$
\limsup_{R\downarrow 0}\frac{1}{|B_R(x_0)|}\int_{B_R(x_0)}|f(x)-f_{B_R(x_0)}|dx=0.
$$
We say that a function $u\in BMO^{-1}(\R^n)$ if there exist $U_j\in
BMO(\R^n)$ such that $u=\sum_{j=1}^n\p_jU_j$. $VMO^{-1}(\R^n)$ is
defined similarly. A remarkable property of $BMO$ function is
\beno
\sup_{R>0,x_0\in\R^n}\frac{1}{|B_R(x_0)|}\int_{B_R(x_0)}|f(x)-f_{B_R(x_0)}|^qdx<\infty.
\eeno
for any $1\le q<\infty$.

Let us conclude this section by recalling the following $\varepsilon$-regularity results.
Here and what follows, we define a solution $u$ to be regular at $z_0=(x_0,t_0)$ if $u\in L^\infty(Q_r(z_0))$ for some $r>0$.

\begin{Proposition}\label{prop:small regularity-CKN}\cite{CKN, LS}
Let $(u,\pi)$ be a suitable weak solution of (\ref{eq:NS}) in $Q_1(z_0)$. There exists an
$\varepsilon_0>0$ such that if
\beno
\int_{Q_1(z_0)}|u(x,t)|^3+|\pi(x,t)|^{3/2}dxdt\leq \varepsilon_0,
\eeno
then $u$ is regular in $Q_{\f 12}(z_0)$. Moreover, $\pi$ can be replaced by
$\pi-(\pi)_{B_r}$ in the integral.
\end{Proposition}

\begin{Proposition}\label{prop:small regularity-GKT}\cite{GKT}
Let $(u,\pi)$ be a suitable weak solution of (\ref{eq:NS}) in $Q_1(z_0)$ and $w=\nabla\times u$.
There exists an $\varepsilon_1>0$ such that if one of the following two conditions holds,
\begin{enumerate}

\item $G(u,p,q;r,z_0)\leq \varepsilon_1$ for any $0<r<\f12$, where $1\leq \frac3p+\frac2q\leq 2$;\vspace{0.1cm}

\item $H(w,p,q;r,z_0)\leq \varepsilon_1$ for any $0<r<\f12$, where $2\leq \frac3{p}+\frac2{q}\leq 3$ and  $(p,q)\neq(1,\infty)$;
\end{enumerate}
then $u$ is  regular at $z_0$.
\end{Proposition}

\section{Interior regularity criterions of suitable weak solution}

The purpose of this section is to establish some interior regularity criterions,
which allow the vertical part of the velocity to be large
under the local scaling invariant norm.
These results improve some classical results and Gustafson-Kang-Tsai's result (Proposition \ref{prop:small regularity-GKT}).
Set $u=(u_h,u_3)$. Let us state our main results.

\begin{Theorem}\label{thm:interior-uball}
Let $(u,\pi)$ be a suitable weak solution of (\ref{eq:NS}) in $Q_1$ and satisfy
\beno
C(u,1)+D(\pi,1)\leq M.
\eeno
Then there exists a positive constant $\varepsilon_2$ depending on $M$ such that if
\beno
C(u_h,1)\leq \varepsilon_2,
\eeno
then $u$ is regular at $(0,0)$.
\end{Theorem}

\begin{Theorem}\label{thm:interior-u}
Let $(u,\pi)$ be a suitable weak solution of (\ref{eq:NS}) in $Q_1$ and satisfy
\beno
G(u,p,q;r)\leq M\quad\textrm{ for any   }0<r<1,
\eeno
where $1\le \frac3p+\frac2q<2$, $1<q\leq \infty$.
There exists a positive constant $\varepsilon_3$ depending on $p, q, M$
such that $(0,0)$ is a regular point if
\beno
G(u_h,p,q;r^*)\leq \varepsilon_3
\eeno
for some $r^*$ with $0<r^*<\min\{\frac1 2, (C(u,1)+D(\pi,1))^{-2}\}$.
\end{Theorem}

\begin{Theorem}\label{thm:interior-gradient}
Let $(u,\pi)$ be a suitable weak solution of (\ref{eq:NS}) in $Q_1$ and satisfy
\beno
H(\nabla u,p,q;r)\leq M\quad \textrm{for any } 0<r<1,
\eeno
where $2\le \frac3{p}+\frac2{q}<3, 1<p\le \infty$.
There exists a positive constant $\varepsilon_4$ depending on $p, q, M$ such that $(0,0)$ is a regular point if
\ben\label{eq:3.6}
H(\nabla u_h,p,q;r^*)\leq \varepsilon_4
\een
for some $r^*$ with $0<r^*<\min\{\frac12, (C(u,1)+D(p,1))^{-2}\}$.
\end{Theorem}

\begin{Remark}
As a special case of Theorem \ref{thm:interior-u}, it follows that
$u$ is regular if
$$|u_3|\le \frac{M}{\sqrt{T-t}},\quad |u_h|\leq\frac{\varepsilon_3}{\sqrt{T-t}},$$
which improves Leray's result \cite{Leray}.
And from Theorem \ref{thm:interior-gradient}, it follows that
$u$ is regular  at $(0,0)$ if for any $0<r<1$,
\beno
r^{-1}\int_{Q_r}|\nabla u_3|^2dxdt\leq M^2,\quad r^{-1}\int_{Q_r}|\nabla u_h|^2dxdt\leq \varepsilon_4^2,
\eeno
which improves Caffarelli-Kohn-Nirenberg's result \cite{CKN}.
\end{Remark}

The proof of Theorem \ref{thm:interior-uball} is based on compactness argument and the following lemma.

\begin{Lemma}\label{lem:c-d}
Let $(u,\pi)$ be a suitable weak solution of (\ref{eq:NS}) in $Q_1$ and $\widetilde{D}(\pi,1)\leq M$.
Then $u$ is regular at $(0,0)$ if
\beno
C(u,r_0)\leq c\varepsilon_0^{9/5}r_0^{8/5}\quad\textrm{ for some }0<r_0\leq 1.
\eeno
Here $c$ is a small constant depending on $M$.
\end{Lemma}

\no{\textbf{Proof.}} By (\ref{eq:pressure}) and H\"{o}lder inequality, for $0<r<{r_0}/4$ we have
\beno
C(u,r)+\widetilde{D}(\pi,r)&\leq&\frac{r_0^2}{r^2}C(u,r_0)+C(\frac{r}{r_0})^{5/2} \widetilde{D}(\pi,r_0)+C\frac{r_0^2}{r^2}C(u,r_0)\\
&\leq& CM\frac{r^{5/2}}{r_0^{9/2}}+C\frac{r_0^2}{r^2}C(u,r_0).
\eeno
Choosing  $r=(\frac{\varepsilon_0}{2CM})^{2/5}r_0^{9/5}$
and  by assumption, we infer that
\beno
C(u,r)+\widetilde{D}(\pi,r)<\varepsilon_0,
\eeno
which implies that $(0,0)$ is a regular point  by Proposition \ref{prop:small regularity-CKN}. \endproof\vspace{0.1cm}

Now let us turn to the proof of  Theorem \ref{thm:interior-uball}.\vspace{0.1cm}

\no{\textbf{Proof of Theorem \ref{thm:interior-uball}}.
Assume that the statement of the proposition is false, then there exist a constant $M$ and a sequence $(u^k,\pi^k)$,
which are suitable weak solutions of (\ref{eq:NS}) in $Q_1$ and singular at $(0,0)$, and satisfies
\beno
C(u^k,1)+D(\pi^k,1)\leq M,\quad C(u_h^k,1)\leq \frac 1k.
\eeno
Then by the local energy inequality, it is easy to get
\beno
A(u^k,3/4)+E(u^k,3/4)\leq C(M),
\eeno
hence by using Lions-Aubin's lemma, there exists a suitable weak solution $(v, \pi')$ of (\ref{eq:NS}) such that (at most up to subsequence),
\beno
u^k\rightarrow v,\quad  u_h^k\rightarrow 0  \quad {\rm in} \quad L^3(Q_\f12),\quad
\pi^k\rightharpoonup \pi'\quad {\rm in} \quad L^\f {3}2(Q_\f12),
\eeno
as $k\rightarrow+\infty$. That is, $v_h=0$, which gives $\partial_3v_3=0$ by $\nabla\cdot v=0$, and hence,
\beno
&\partial_tv_3+\p_3 \pi'-\triangle v_3=0,\quad -\triangle\pi'=0,\quad \textrm{or}\\
&\p_tv+\na \pi'-\triangle v=0,
\eeno
which implies that $|v|\leq C(M)$ in $Q_{1/4}$ by the classical result of linear Stokes equation(see \cite{chen2} for example).
However, $(0,0)$ is a singular point of $u^k$, hence by Lemma \ref{lem:c-d}, for any $0<r<1/4$,
\beno
c\varepsilon_0^{9/5}r^{8/5}&\leq& \lim_{k\rightarrow\infty}r^{-2}\int_{Q_r}|u^k|^3dxdt\\
&\leq& \lim_{k\rightarrow\infty} C(v,r)
\leq C(M)r^{3},
\eeno
which is a contradiction by letting $r\rightarrow 0$.\endproof

\vspace{0.2cm}

The proof of Theorem \ref{thm:interior-u} is motivated by \cite{Se3} and based on the blow-up argument.
\vspace{0.1cm}

\no{\textbf{Proof of Theorem \ref{thm:interior-u}}}.
Assume that the statement of the proposition is false, then there exist constants $p, q, M$ and a sequence $(u^k,\pi^k)$,
which are suitable weak solutions of (\ref{eq:NS}) in $Q_1$ and singular at $(0,0)$, and satisfy
\beno
&&G(u^k,p,q;r)\leq M \quad {\rm for \,\,all}\quad 0<r<1,\\
&&G(u_h^k,p,q;r_k)\leq \frac1k,
\eeno
where $0<r_k<\min\{\frac12, (C(u^k,1)+D(\pi^k,1))^{-2}\}$. Then it follows from Lemma \ref{lem:invariant} that
\beno
A(u^k,r)+E(u^k,r)+D(\pi^k,r)\leq C(M,p,q)
\eeno
for any $0<r<r_k.$

Set $v^k(x,t)=r_ku^k(r_kx,r_k^2t), q^k(x,t)=r_k^2\pi^k(r_kx,r_k^2t)$. Then
\beno
A(v^k,r)+E(v^k,r)+D(q^k,r)\leq C(M,p,q)
\eeno
for any $0<r<1$. Lions-Aubin's lemma ensures that there exists a suitable weak solution $(\bar{v}, \bar{\pi})$ of (\ref{eq:NS})
such that (at most up to subsequence),
\beno
&&v^k\rightarrow \bar{v}\quad {\rm in} \quad L^3(Q_\f12),\quad q^k\rightharpoonup \bar{q}\quad {\rm in} \quad L^\f {3}2(Q_\f12),\\
&& v_h^k\rightharpoonup 0\quad {\rm in} \quad L^q((-\f14,0); L^p(B_\f12)),
\eeno
as $k\rightarrow+\infty$.
Then we have $\bar{v}_h=0$ and
$$\partial_t\bar{v}_3+\p_3 \bar{q}-\triangle \bar{v}_3=0,$$
which implies that $|\bar{v}_3|\leq C(M)$ in $Q_{\f14}$.
However, $(0,0)$ is a singular point of $v^k$,
hence by Proposition \ref{prop:small regularity-CKN} and (\ref{eq:pressure1}), for any $0<r<1/4$,
\beno
\varepsilon_0 &\leq&\liminf_{k\rightarrow\infty}r^{-2}\int_{Q_r}|v^k|^3+|q^k|^{3/2}dxdt\\
&\leq& C\liminf_{k\rightarrow\infty}\Big(C(\bar{v},r)+\frac{r}{\rho}D(q^k,\rho)+(\frac{\rho}{r})^2C(v^k,\rho)\Big)\\
&\leq& C\Big(C(\bar{v},r)+\frac{r}{\rho}+(\frac{\rho}{r})^2C(\bar{v},\rho)\Big)\\
&\leq& Cr^{1/2}\quad (\textrm{by choosing}\quad \rho=r^{\f12}),
\eeno
which is a contradiction if we take $r$ small enough.\endproof

\no{\bf Proof of Theorem \ref{thm:interior-gradient}}. Without loss of generality, let us assume that
\beno
\f 83<\f 3p+\f 2q<3.
\eeno
The other case can be reduced to it by H\"{o}lder inequality. By Lemma \ref{lem:invariant}, we have
\beno
A(u,r)+E(u,r)+D(\pi,r)\leq C(M)\big(r^{1/2}\big(C(u,1)+D(\pi,1)\big)+1\big)\leq C(M),
\eeno
for any $0<r\leq r_1\triangleq\min\{\frac12, \big(C(u,1)+D(p,1)\big)^{-2}\}$. This together with interpolation inequality gives
\ben\label{eq:3.7}
C(u,r)\leq C(M)\quad \textrm{for any  } 0<r\leq r_1.
\een
We get by Poinc\'{a}re inequality that
\beno
\widetilde{G}(u_h,p_1,q_1;r)\le C\widetilde{H}(\na u_h,p,q;r),
\eeno
where $p_1=\f {3p} {3-p}, q_1=q$, hence it follows from (\ref{eq:7.9}) and (\ref{eq:3.6}) that
\beno
\widetilde{C}(u_h,r^*)&\leq& C(M)\big(A(u_h,r^*)+E(u_h,r^*)\big)^{\frac{1-3\delta}{1-2\delta}}
\widetilde{G}(u_h,p_1,q_1;r^*)^{\frac{1}{1-2\delta}},\\
&\leq& C(M)\widetilde{G}(u_h,p_1,q_1;r^*)^{\frac{1}{1-2\delta}}\le C(M)\varepsilon_4^{\frac{1}{1-2\delta}},
\eeno
where $\delta=2-\frac3{p_1}-\frac2{q_1}\in (0,\f13)$, hence by (\ref{eq:3.7}) for $0<r<r^*$
\beno
C(u_h,r)\leq C(\frac{r}{r^*})C(u_h,r^*)+C(\frac{r^*}{r})^2\widetilde{C}(u_h,r^*)\leq C(M)
\big(\frac{r}{r^*}+(\frac{r^*}{r})^2\varepsilon_4^{\frac{1}{1-2\delta}}\big).
\eeno
Taking $r$ small enough, and then $\varepsilon_4$ small enough such that
\beno
C(u_h,r)\leq \varepsilon_3^3.
\eeno
Then the result follows from Theorem \ref{thm:interior-u}.\endproof

\section{Applications of interior regularity criterions}

\subsection{Ladyzhenskaya-Prodi-Serrin's criterion}

Using the interior regularity criterions established in Section 3, we present Ladyzhenskaya-Prodi-Serrin's type criterions
in Lorentz spaces.

\begin{Theorem}\label{thm:serrin}
Let $u$ be a Leray-Hopf weak solution of (\ref{eq:NS}) in $\R^3\times(-1,0)$. Assume that $u$ satisfies
\ben\label{eq:4.8}
 \|u\|_{L^{q,\infty}((-1,0); L^{p,\infty}(\R^3))}\leq M,\quad \|u_h\|_{L^{q,\ell}((-1,0); L^{p,\infty}(\R^3))}<\infty,
\een
where $\frac3p+\frac2q=1,$ $ 3<p\leq \infty$, and $1\leq \ell<\infty$. Then $u$ is regular in $\R^3\times(-1,0]$.
For $\ell=\infty$ or $p=3$, the same result holds if the second condition of (\ref{eq:4.8}) is replaced by
\beno
\|u_h\|_{L^{q,\infty}((-1,0); L^{p,\infty}(\R^3))}\leq \varepsilon_5,
\eeno
where $\varepsilon_5$ is a small constant depending on $M$.
\end{Theorem}

\begin{Remark}
For $\ell=\infty$, we improve Kim-Kozono's result \cite{Kim} and He-Wang's result \cite{He}, where the smallness of all components of the velocity
is imposed. In general case, we improve Sohr's result \cite{Sohr} by allowing the vertical part of the velocity to fall in weak $L^p$ space.
\end{Remark}

\begin{Remark}
Under the condition (\ref{eq:4.8}), it can be verified that Leray-Hopf weak solution is suitable weak solution.
We left it to the interested readers.
\end{Remark}

The proof is based on the following lemma.

\begin{Lemma}\label{lem:local bound}
Assume that $u$ satisfies
\beno
 \|u\|_{L^{q,\infty}((-1,0); L^{p,\infty}(\R^3))}\leq m,
\eeno
where $\frac3p+\frac2q=1,$ $ 3\leq p\leq \infty$. Then for any $0<r<1$ and $0<\epsilon<1$, there hold
\beno
&&G(u,\frac{9}{10}p,\frac{4}{5}q;r)\leq C\epsilon^{\frac{4q}{5}}+C\epsilon^{-\frac{q}{5}}m^{q},\quad 3<p<\infty,  \\
&&A(u,r)\leq C\epsilon^2+C\epsilon^{-1}m^3,\quad p=3, \\
&&G(u,\infty,\f32;r)\leq C\epsilon^{3/2}+C\epsilon^{-1/2}m^2,\quad p=\infty.
\eeno
\end{Lemma}

\no{\bf Proof}. First we consider the case of $3<p<\infty$. Using
the definition of Lorentz space, we infer that
\beno
&&r^{(\frac45-\f 8{3p})q-2}\int_{-r^2}^0\Big(\int_{B_r}|u|^{\f {9}{10}p}dx\Big)^{\f {8q} {9p}}dt\\
&&\leq Cr^{(\frac45-\f 8{3p})q-2}\int_{-r^2}^0\Big(\int_0^{\infty}\sigma^{\f {9} {10}p-1}|\{x\in B_r; |u(x,t)|>\sigma\}| d\sigma \Big)^{\f {8q} {9p}}dt\\
&&\leq Cr^{(\frac45-\f 8{3p})q-2}\Big(r^2 (r^3R^{\f {9}{10}p})^{\f {8q} {9p}}
+\int_{-r^2}^0\Big(\int_R^{\infty}\sigma^{\f {9} {10}p-1}|\{x\in B_r; |u(x,t)|>\sigma\}| d\sigma \Big)^{\f {8q} {9p}}dt\Big)\\
&&\leq Cr^{(\frac45-\f 8{3p})q-2}\Big(r^2 (r^3R^{\f {9}{10}p})^{\f {8q}
{9p}}+R^{-\f {4q}{45}}\int_{-r^2}^0\|u(\cdot,t)\|_{L^{p,\infty}}^{\f {8q} {9}}dt\Big)\\
&&\leq Cr^{(\frac45-\f 8{3p})q-2}\Big(r^2 (r^3R^{\f {9}{10}p})^{\f {8q}
{9p}}+R^{-\f {4q}{45}}r^{-(1-\frac3p)\f {8q} 9+2}I(r)\Big)\\
&&\leq C\epsilon^{\f {4q} 5}+C\epsilon^{-\f
{4q}{45}}I(r),
\eeno
where we take $R=\epsilon r^{-1}$ and the estimate of $I(r)$ is given by
 \beno
&&I(r)\equiv r^{(1-\frac3p)\f {8q} 9-2}\int_{-r^2}^0\|u(\cdot,t)\|_{L^{p,\infty}(B_1)}^{\f {8q} 9}dt\\
&&\leq Cr^{(1-\frac3p)\f {8q} 9-2}\int_{0}^{\infty}\sigma^{\f {8q} 9-1}|\{t\in(-r^2,0); \|u(\cdot,t)\|_{L^{p,\infty}}>\sigma \}|d\sigma\\
&&\leq Cr^{(1-\frac3p)\f {8q} 9-2}\Big(R^{\f {8q} 9}r^2+\int_{R}^{\infty}\sigma^{\f {8q} 9-1}|\{t\in(-r^2,0); \|u(\cdot,t)\|_{L^{p,\infty}}>\sigma \}|d\sigma\Big)\\
&&\leq Cr^{(1-\frac3p)\f {8q} 9-2}\Big(R^{\f {8q} 9}r^2+R^{-\f q9}\|u\|_{L^{q,\infty}(-1,0; L^{p,\infty}(B_1))}^q\Big)\\
&&\leq Cr^{(1-\frac3p)\f {8q} 9-2}\Big(R^{\f {8q} 9}r^2+R^{-\f q9}m^q\Big)\\
&&\leq C\epsilon^{\f {8q} 9}+C\epsilon^{-\f {q} 9}m^q\quad
(R=\epsilon r^{\frac3p-1}).
\eeno
This gives the first inequality. For $p=3$, we consider
\beno
\sup_{-r^2<t<0}r^{-1}\int_{B_r}|u|^2dx
&\leq& C\sup_{-r^2<t<0}r^{-1}\int_0^{\infty}\sigma|\{x\in B_r; |u(x,t)|>\sigma\}|d\sigma\\
&\leq& C\sup_{-r^2<t<0}r^{-1}\Big(R^2r^3+\int_R^{\infty}\sigma|\{x\in B_r; |u(x,t)|>\sigma\}|d\sigma\Big)\\
&\leq& C\sup_{-r^2<t<0}r^{-1}\Big(R^2r^3+R^{-1}\|u(\cdot,t)\|_{L^{3,\infty}}^3\Big),
\eeno
which gives the second inequality by taking $R=\epsilon r$. Let $g(t)=\|u(\cdot,t)\|_{L^\infty(B_1)}$. Then we have
\beno
r^{-1/2}\int_{-r^2}^0g(t)^{3/2}dt
&\leq& Cr^{-1/2}\int_0^{\infty}\sigma^{\f12}|\{t\in (-r^2,0); |g(t)|>\sigma\}|d\sigma\\
&\leq& Cr^{-1/2}\Big( R^{\f32}r^2+\int_R^{\infty}\sigma^{\f12}|\{t\in (-r^2,0); |g(t)|>\sigma\}|d\sigma\Big)\\
&\leq& Cr^{-\f 12}\Big(R^{\f32}r^2+R^{-\f12}m^2\Big),
\eeno
which gives the third inequality by taking $R=\epsilon r$.
\endproof

\vspace{0.2cm}

\no{\bf Proof of Theorem \ref{thm:serrin}}. By translation invariance and Theorem \ref{thm:interior-u},
it suffices to show that
\ben\label{eq:5.3}
G(u,p_1,q_1;r)\leq M,\quad G(u_h,p_1,q_1;r)\leq \varepsilon_3,
\een
for any $0<r<1/2$ and some $(p_1,q_1)$ with $1\leq \frac3{p_1}+\frac2{q_1}<2$.
For $3<p<\infty$, let $p_1=\frac{9}{10}p$ and $q_1=\frac{4}{5}q$, then $\frac3{p_1}+\frac2{q_1}<\frac54<2$.
For $\ell<\infty$, we have  $\| u_h\|_{L^{q,\ell}(-r^2,0; L^{p,\infty}(B_r))}\rightarrow 0$ as $r\rightarrow 0$.
Hence by Lemma \ref{lem:local bound}, the condition (\ref{eq:5.3}) holds
if we take $\epsilon$ small enough, and then take $r$ small enough.
The proof of the other cases is similar. We omit the details.
\endproof

\subsection{Escauriza-Seregin-\v{S}ver\'{a}k's criterion}

The following theorem improves Escauriza-Seregin-\v{S}ver\'{a}k's criterion by noting the inclusion
\beno
L^3(\R^3)\subset L^{3,\ell}(\R^3)\quad\textrm{ for }\ell>3 \quad \textrm{and} \quad   L^3(\R^3)\subset VMO^{-1}(\R^3).
\eeno

\begin{Theorem}\label{thm:main2}
Let $(u,\pi)$ be a suitable weak solution of (\ref{eq:NS}) in $\R^3\times(-1,0)$. If
\beno
\|u_h\|_{L^{\infty}((-1,0); L^{3,\ell}(\R^3))}+\|u_3\|_{L^{\infty}((-1,0); BMO^{-1}(\R^3))}=M<\infty,
\eeno
for some $\ell<\infty$, and  $u_3(x,t)\in VMO^{-1}(\R^3)$ for $t\in (-1,0]$,
then $u$ is regular in $\R^3\times(-1,0]$.
\end{Theorem}

We need the following lemma, which gives a bound of local scaling invariant energy.

\begin{Lemma}\label{lem:energy}
Under the assumptions of Theorem \ref{thm:main2}, there holds
 \beno
A(u,r)+E(u,r)+D(\pi,r)\leq C(M, C(u,1), D(\pi,1))\quad\textrm{ for any  } \,0<r<1/2.
\eeno
\end{Lemma}

\no{\bf Proof.}
Let $\zeta(x,t)$ be a smooth function with $\zeta\equiv 1$ in $Q_r$
and $\zeta=0$ in $Q_{2r}^c$. Since $u_3\in
L^{\infty}(-1,0;BMO^{-1}(\R^3))$, there exists $U(x,t)\in
L^{\infty}(-1,0;BMO(\R^3))$ such that $u_3=\na\cdot U$. We have by
H\"{o}lder inequality that \beno
&&r^{-2}\int_{Q_{2r}}|u_3|^3\zeta^2dxdt\\
&&=r^{-2}\int_{Q_{2r}}\sum_{j=1}^3\partial_j U_j\cdot u_3 |u_3|\zeta^2dxdt\\
&&\leq 6r^{-2}\int_{Q_{2r}}|U-U_{B_{2r}}|(|\nabla u_3| |u_3|+|u_3|^2|\nabla\zeta|)dxdt\\
&&\leq 6r^{-2}\Big(\int_{Q_{2r}}|U-U_{B_{2r}}|^6dxdt\Big)^{1/6}\Big(\int_{Q_{2r}}|\nabla u_3|^2dxdt\Big)^{1/2}\Big(\int_{Q_{2r}}|u_3|^3dxdt\Big)^{1/3}\\
&&\quad+12r^{-3}\Big(\int_{Q_{2r}}|U-U_{B_{2r}}|^3dxdt\Big)^{1/3}\Big(\int_{Q_{2r}}|u_3|^3dxdt\Big)^{2/3},
\eeno
which implies that
\ben\label{eq:4.10}
C(u_3,r)\leq
C(M)\big(E(u,2r)^{1/2}C(u,2r)^{1/3}+ C(u,2r)^{2/3}\big).
\een
On the other hand, we have by Lemma \ref{lem:local bound} that
\beno
A(u_h,r)\leq C(M)\quad \textrm{for any } 0<r<1,
\eeno
which along with the interpolation
inequality gives
\ben\label{eq:4.11}
C(u_h,r)\leq A(u,r)^{3/4}\big(E(u,r)+A(u,r)\big)^{3/4}\leq
C(M)\big(E(u,r)+A(u,r)\big)^{3/4}.
\een
We infer from (\ref{eq:4.10}) and (\ref{eq:4.11}) that
\beno
C(u,r)\leq C(M)\big(E(u,2r)+A(u,2r)+C(u,2r)\big)^{5/6}.
\eeno
With this, following the proof of Lemma \ref{lem:invariant}, we conclude the result. \endproof
\vspace{0.2cm}

\no{\bf Proof of Theorem \ref{thm:main2}.}\,
Following \cite{ESS}, the proof is based on the blow-up analysis and unique continuation theorem.
Without loss of generality, assume that $(0,0)$ is a singular point.
Then by Theorem \ref{thm:interior-u}, there exists a sequence of $r_k\downarrow 0$ such that
\ben\label{eq:4.2}
r_k^{-2}\int_{Q_{r_k}}|u_h|^3dxdt\geq \varepsilon_1.
\een
Let $u^k(x,t)=r_ku(r_kx,r_k^2t)$ and $\pi^k(x,t)=r_k^2\pi(r_kx,r_k^2t)$. Then
for any $a>0$ and $k$ large enough, it follows from Lemma \ref{lem:energy} that
\beno
A(u^k,a)+E(u^k,a)+C(u^k,a)+D(\pi^k,a)\leq C(M,D(\pi,1)).
\eeno
Using Lions-Aubin lemma, there exists $(v,\pi')$ such that for any $a,T>0$ (up to subsequence)
\beno
&&u^k\rightarrow v \quad {\rm in} \quad L^3(B_{a}\times (-T,0)),\\
&&u^k\rightarrow v \quad {\rm in} \quad C([-T,0];L^{9/8}(B_{a})),\\
&&\pi^k\rightharpoonup \pi' \quad {\rm in} \quad L^\f 32(-T,0;L^\infty(B_a)),
\eeno
as $k\rightarrow+\infty$ (see the proof of Theorem 4.1 in \cite{WZ} for the details). Furthermore, there hold
\ben\label{eq:4.15}
\|v_h\|_{L^{\infty}(-a^2,0; L^{3,\ell}(\R^3))}\leq \sup_{k}\|u_h^k\|_{L^{\infty}(-a^2,0; L^{3,\ell}(\R^3))}\leq M,
\een
and for any $z_0=(x_0,t_0)\in (-T+1,0)\times \R^3$,
\ben\label{eq:4.16}
A(v,1;z_0)+E(v,1;z_0)+C(v,1;z_0)+D(\pi',1;z_0)\leq C(M,D(p,1)).
\een
Due to (\ref{eq:4.15}) and (\ref{eq:lorentz}),  we infer that
\beno
\int_{Q_1(z_0)}|v_h|^2dxdt\rightarrow0, \quad {\rm as}\,\, z_0\rightarrow\infty,
\eeno
which  along with (\ref{eq:4.16}) implies that
\beno
\int_{Q_1(z_0)}|v_h|^3dxdt\rightarrow0, \quad {\rm as}\,\, z_0\rightarrow\infty.
\eeno
Hence by Theorem \ref{thm:interior-uball}, there exists $R>0$ such that
\beno
|v(x,t)|+|\nabla v(x,t)|\leq C, \quad (t,x)\in (-T+1,0)\times \R^3\backslash{B_R}.
\eeno

Due to $u_h(x,0)\in L^{3,\ell}$, we infer that
\beno
\int_{B_a}|v_h(x,0)| dx
&\leq& \int_{B_a}|v_h(x,0)-u_h^k(x,0)|dx+\int_{B_a}|u_h^k(x,0)| dx\\
&\leq&\int_{B_a}|v_h(x,0)-u_h^k(x,0)|dx+r_k^{-2}\int_{B_{ar_k}}|u_h(y,0)|dy\\
&\le&\int_{B_a}|v_h(x,0)-u_h^k(x,0)|dx+C\|u_h(0)\|_{L^{3,\ell}(B_{ar_k})}\\
&\longrightarrow& 0,\quad\textrm{ as } k\rightarrow\infty,
\eeno
which implies $v_h(x,0)=0$ a.e. $\R^3$. And due to $u_3(x,0)\in VMO^{-1}(\R^3)$, we have $v_3(x,0)=0$
(see Theorem 4.1 in \cite{WZ}).

Let $w=\nabla\times v$, then $w(x,0)=0$ and
$$|\partial_t w-\triangle w|\leq C(|w|+|\nabla w|),\quad (-T+1,0)\times \R^3\backslash{B_R}. $$
By the backward uniqueness property of parabolic operator \cite{ESS}, we have
$w=0$ in $(-T+1,0)\times \R^3\backslash{B_R}.$ Similar arguments as in \cite{ESS}, using spacial unique continuation
we have $w\equiv 0$ in $(-T+1,0)\times \R^3$, which implies $\triangle v\equiv 0$ in $(-T+1,0)\times \R^3$,
hence $v_h\equiv 0$ in $(-T+1,0)\times \R^3$,
since $v_h(\cdot,t)\in L^{3,\ell}$. This is a contradiction to (\ref{eq:4.2}).\endproof

\section{The number of singular points}

\subsection{An improved $\varepsilon$-regularity criterion}

We need the following improved version, which may be independent of interest.

\begin{Proposition}\label{prop:small regularity-new}
Let $(u,\pi)$ be a suitable weak solution of (\ref{eq:NS}) in $Q_1(z_0)$.
There exists an $\varepsilon_6>0$ such that if
\beno
&&\sup_{t\in [-1+t_0,t_0]}\int_{B_1(x_0)}|u(x,t)|^2dx+\int_{-1+t_0}^{t_0}\big(\int_{B_1(x_0)}|u(x,t)|^4dx\big)^{\f12}dt\\
&&\quad+\int_{-1+t_0}^{t_0}\big(\int_{B_1(x_0)}|\pi(x,t)|^2dx\big)^{\f12}dt\leq \varepsilon_6,
\eeno
then $u$ is regular in $Q_{\f1 2}(z_0)$.
\end{Proposition}

\begin{Remark}
Due to Lemma \ref{lem:pressure}, the above norm of the pressure can be replaced by $L^1(Q_1(z_0))$ norm.
A slightly different version of Proposition \ref{prop:small regularity-new} was obtained by Vasseur \cite{Vasseur},
who used the De Giorgi iterative method.
\end{Remark}

\no{\bf Proof.}\,By Proposition \ref{prop:small regularity-GKT} and translation invariance, it suffices to prove that
\ben\label{eq:regular}
A(u,r)+E(u,r)\leq \varepsilon_6^\f12\le \varepsilon_1^2
\een
for any $0<r<1/2$. Set $r_n=2^{-n}$, where $n=1,2,\cdots$.
First of all, (\ref{eq:regular}) holds for $r=r_1$ by local energy inequality.
Suppose that (\ref{eq:regular}) holds for $r_k$ with $k\leq n-1$. We need to show that
\ben\label{eq:regular2}
A(u,r_n)+E(u,r_n)\leq \varepsilon_6^\f12.
\een
Let $\phi_n=\chi\psi_n$, where $\chi$ is a cutoff function which equals $1$ in $Q_{1/4}$ and vanishes outside of $Q_{1/3}$,
and $\psi_n$ is as follows:
\beno
\psi_n=(r_n^2-t)^{-3/2}e^{-\frac{|x|^2}{4(r_n^2-t)}}.
\eeno
Direct computations show that $\phi_n\geq0$ and
$$(\partial_t+\triangle)\phi_n=0\quad \textrm{in}\quad  Q_{1/4},$$
$$|(\partial_t+\triangle)\phi_n|\leq C_1\quad \textrm{in}   \quad  Q_{1/3},$$
$$C_1^{-1}r_n^{-3}\leq \phi_n\leq C_1 r_n^{-3},\quad |\nabla\phi_n|\leq C_1r_n^{-4}\quad \textrm{on}\quad  Q_{r_n} \quad n\geq 2,$$
$$\phi_n\leq C_1r_k^{-3},\quad |\nabla\phi_n|\leq C_1r_k^{-4}\quad \textrm{on}\quad  Q_{r_{k-1}}/{Q_{r_k}} \quad 1<k\leq n.$$
Using $\phi_n$ as a test function in the local energy inequality, we get
\beno
&&\sup_{-r_n^2<t<0}r_n^{-1}\int_{B_{r_n}}|u(x,t)|^2 dx+r_n^{-1}\int_{Q_{r_n}}|\nabla u|^2dxdt\\
&&\leq C_1^2r_n^2\int_{Q_1}|u|^2dxdt+C_1r_n^2\int_{Q_1}|u|^3|\nabla\phi_n|dxdt+C_1r_n^2\big|\int_{Q_1}\pi(u\cdot\nabla\phi_n)dxdt\big|\\
&&\eqdef I_1+I_2+I_3.
\eeno

Firstly, we have by assumption that
\beno
I_1\le C_1^2 r_n^2\varepsilon_6.
\eeno
Recall that the following well-known interpolation inequality from \cite{CKN}: for $\rho\ge r>0$
\beno
 C(u,r)\leq C(\frac{\rho}{r})^3A(u,\rho)^{3/4}E(u,\rho)^{3/4}+C(\frac{r}{\rho})^{3}A(u,\rho)^{3/2},
\eeno
from which and the induction assumption, it follows that
\beno
I_2 &\leq& C_1^2r_n^2\sum_{k=1}^nr_k^{-4}\int_{Q_{r_k}}|u|^3dxdt\\
&\leq & Cr_n^2\sum_{k=1}^nr_k^{-2}\varepsilon_6^{3/4}
\leq  C\varepsilon_6^{3/4}.
\eeno
To estimate $I_3$, we choose $\chi_k$ to be a cutoff function, which vanishes
outside of $Q_{r_k}$ and equals 1 in $Q_{7/8r_k}$, and $|\nabla\chi_k|\leq Cr_k^{-1}$.
We have by the induction assumption that
\beno
I_3 &\leq & C_1r_n^2\sum_{k=1}^{n-1}\big|\int_{Q_{1}}\pi(u\cdot\nabla((\chi_k-\chi_{k+1})\phi_n))dxdt\big|
+C_1r_n^2\big|\int_{Q_{1}}\pi u\cdot\nabla(\chi_n\phi_n)dxdt\big|\\
&\leq & C_1r_n^2\sum_{k=1}^{n-1}\big|\int_{Q_{1}}(\pi-(\pi)_{B_k})u\cdot\nabla((\chi_k-\chi_{k+1})\phi_n)dxdt\big|\\
&&+C_1r_n^2\big|\int_{Q_{1}}(\pi-(\pi)_{B_n})u\cdot\nabla(\chi_n\phi_n)dxdt\big|\\
&\leq& Cr_n^2\sum_{k=3}^{n}r_k^{-4}\int_{Q_{r_k}}|(\pi-(\pi)_{B_k})u|dxdt+Cr_n^2\int_{Q_1}|u||\pi|dxdt\\
&\leq& Cr_n^2\sum_{k=3}^{n}r_k^{-2}\varepsilon_6^{1/4}\widetilde{H}(\pi,2,1;r_k)+C\varepsilon_6^{3/2},
\eeno
and by Lemma \ref{lem:pressure} and interpolation inequality,  we get
\beno
\widetilde{H}(\pi,2,1;\theta^j)&\leq& C\theta\widetilde{H}(\pi,1,1;\theta^{j-1})+C\theta^{-\frac32}G(u,4,2;\theta^{j-1})^2\\
&\leq& (C\theta)^j\widetilde{H}(\pi,1,1;1)+C\theta^{-\frac32}\sum_{\ell=1}^j(C\theta)^{\ell-1}G(u,4,2;\theta^{j-\ell})^2\\
&\leq &(C\theta)^j\varepsilon_6+C\theta^{-\frac32}\sum_{l=1}^j(C\theta)^{\ell-1}\varepsilon_6^\f12\\
&\leq & C\varepsilon_6^\f12,
\eeno
where we take $\theta$ such that $C\theta<\f12$ and $j$ satisfies $\theta^j\geq r_{n}$.
This gives
\beno
I_3\leq  C\varepsilon_6^{3/4}.
\eeno
Summing up the estimates for $I_1-I_3$ and taking $\varepsilon_6$ small enough, we conclude (\ref{eq:regular2}).
\endproof

\subsection{The number of singular points}

\begin{Theorem}\label{thm:singular point}
Let $u$ be a Leray-Hopf weak solution in $\R^3\times(-1,0)$ and satisfy
\ben\label{eq:6.3}
\|u\|_{L^{q,\infty}(-1,0; L^{p}(\R^3))}=M<\infty,
\een
where $\frac3p+\frac2q=1$, $3<p<\infty$.
Then the number of singular points of $u$ is finite at any time $t\in (-1,0]$, and the number depends on $M$.
\end{Theorem}

\begin{Remark}
The case of $(p,q)=(3,\infty)$ has been proved by Neustupa \cite{Ne} and Seregin \cite{Seregin-CPAM}.
In fact, the solution is regular in this case \cite{ESS}.
A special case satisfying (\ref{eq:6.3}) is
\beno
\|u(t)\|_{L^p(\R^3)}\leq M(-t)^{\frac{3-p}{2p}}.
\eeno
Note that the solution is regular if $M$ is small, which was proved by Leray \cite{Leray}.
\end{Remark}

\begin{Lemma}\label{lem:regular2}
Let $(u,\pi)$ be a suitable weak solution of (\ref{eq:NS}) in $Q_1$ and satisfy
\ben\label{eq:6.4}
\|u\|_{L^{q,\infty}(-1,0; L^{p}(B_1))}<M,
\een
where $\frac3p+\frac2q=1$ and $3<p<\infty$.
There exists $\varepsilon_7>0$ depending on $M, C(u,1), D(\pi,1)$ such that $u$ is regular at $(0,0)$ if
\ben\label{eq:6.5}
\|u\|^2_{L_t^{q_0}L_x^p}(Q_1)+\|\pi\|_{L_t^{q_0/2}L_x^{p/2}}(Q_1)\leq \varepsilon_7.
\een
where $q_0=3$ for $3<p<9$ and $q_0=\frac{q+2}{2}$ for $p\geq 9$.
\end{Lemma}

\no{\bf Proof.}
For $p\in (3,9)$, the result follows from Proposition \ref{prop:small regularity-CKN} and (\ref{eq:6.5}).
Now we assume $p\geq 9$.
Similar to the proof of Lemma \ref{lem:local bound}, we can infer from (\ref{eq:6.4}) that
\ben\label{eq:6.6}
G(u,p,q_0;r)\leq C(M) \quad\textrm{ for any }0<r<1,
\een
which along with Lemma \ref{lem:invariant} gives
\ben\label{eq:16}
A(u,r)+E(u,r)+D(\pi,r)\leq C(M, C(u,1),D(\pi,1))\quad\textrm{ for any  } 0<r<1/2.
\een
Due to $p\ge 9$, hence $q_0>2$, H\"{o}lder inequality gives
\ben\label{eq:5.20}
G(u,4,2;r)\le CG(u,p,q_0;r).
\een

Let $\zeta$ be a cutoff function, which vanishes
outside of $Q_{\rho}$ and equals 1 in $Q_{\rho/2}$, and satisfies
$$|\nabla\zeta|\leq C_1\rho^{-1},\quad |\partial_t\zeta|,|\triangle\zeta|\leq C_1\rho^{-2}.$$
Define the backward heat kernel as
$$\Gamma(x,t)=\frac{1}{4\pi(r^2-t)^{3/2}}e^{-\frac{|x|^2}{4(r^2-t)}}.$$
Taking the test function $\phi=\Gamma\zeta$ in the local energy inequality, and noting $(\partial_t+\triangle)\Gamma=0$,
we obtain
\beno
\sup_t\int_{B_{\rho}}|u|^2\phi dx+\int_{Q_{\rho}}|\nabla u|^2\phi dxdt
\leq \int_{Q_{\rho}}\big(|u|^2(\triangle\phi+\partial_t\phi)+u\cdot\nabla\phi(|u|^2+\pi)\big)dxdt.
\eeno
This implies that
\beno
A(u,r)+E(u,r)\leq C(\frac{r}{\rho})^2\Big(\rho^{-3}\int_{Q_{\rho}}|u|^2dxdt+C(u,\rho)+
\rho^{-2}|u||\pi-(\pi)_{B_{\rho}}|dxdt\Big).
\eeno
While by (\ref{eq:6.6}) and (\ref{eq:5.20}), we have
$$C(u,\rho)\leq A(u,\rho)^{1/2}G(u,4,2;\rho)^2\leq C(M)A(u,\rho)^{1/2}.$$
And we get by Lemma 6.1 that
\beno
\widetilde{H}(\pi,2,1;r)&\leq& C(\frac{\rho}{r})^{\frac{3}{2}} G(u,4,2;\rho)^2+C (\frac{r}{\rho})\widetilde{H}(\pi,1,1;\rho)\\
&\leq& C(M)(\frac{\rho}{r})^{\frac{3}{2}} +C(\frac{r}{\rho})\widetilde{H}(\pi,1,1;\rho),
\eeno
which gives by a standard iteration that
\beno
\tilde{H}(\pi,2,1;r)\leq C(M)\quad\textrm{for } 0<r<1/2.
\eeno
Hence, we have
$$
\rho^{-2}\int_{Q_{\rho}}|u||\pi-\pi_{B_{\rho}}|dxdt\leq C A(u,\rho)^{1/2}\tilde{H}(\pi,2,1;\rho)
\le C(M)A(u,\rho)^{1/2}.
$$
Let $F(r)=A(u,r)+E(u,r)+\tilde{H}(\pi,2,1;r)^2$. Then we conclude
\ben\label{eq:17}
F(r)\leq C (\frac{r}{\rho})^2F(\rho)+C(M)(\frac{r}{\rho})^2+C(\frac{\rho}{r})^{3} G(u,4,2;\rho)^4.
\een
Letting $\rho=1$, and taking $r$ small and then $\varepsilon_7$ small,
we infer from (\ref{eq:17}), (\ref{eq:16}) and Proposition \ref{prop:small regularity-new} that $(0,0)$ is a regular point.
The proof is completed.\endproof
\vspace{0.2cm}

Now we are in position to prove Theorem \ref{thm:singular point}.\vspace{0.1cm}

\no\textbf{Proof of Theorem \ref{thm:singular point}}.
We denote $z_1=(x_1, t_0),\cdots,z_K=(x_K,t_0)$ by the singular points of the solution at $t=t_0$.
Then Lemma \ref{lem:regular2} implies that at every singular point we have
\beno
G(u,p,q_0;r,z_l)^2+H(\pi,p/2,q_0/2;r,z_l)> \varepsilon_7,\quad \textrm{for any }0<r<1,
\eeno
where $l=1,\cdots,K$. We choose $r_0>0$ small such that $B_r(x_i)\cap B_r(x_j)=\emptyset$ for $i\neq j$ and all $0<r\leq r_0$.
Taking $r=\theta^kr_0$ and $\rho=\theta^{k-1}r_0$ in (\ref{eq:pressure1}), we find
\beno
&&H(\pi,p/2,q_0/2;\theta^kr_0)^{q_0/2}\\
&&\leq C\theta^{q_0-2}H(\pi,p/2,q_0/2;\theta^{k-1}r_0)^{q_0/2}+C\theta^{-2-\frac{3q_0}p+q_0} G(u,p,q_0;\theta^{k-1}r_0)^{q_0}\\
&&\leq (C\theta^{q_0-2})^kH(\pi,p/2,q_0/2;r_0)^{q_0/2}
+C\theta^{-2-\frac{3q_0}p+{q_0}}\sum_{i=0}^{k-1}(C\theta^{q_0-2})^{k-i-1}G(u,p,q_0;\theta^{i}r_0)^{q_0}.
\eeno
Now, for $\frac{q_0}{p}\geq 1$, noting that
$$
\sum_{l=1}^Ka_l^{\f {q_0} p}\le \big(\sum_{l=1}^Ka_l\big)^{\f {q_0} p},\quad a_l\ge 0,
$$
we deduce by (\ref{eq:lorentz}) with $R=r^{-2/q}\|u\|_{L^{q,\infty}(-1,0;L^p(\R^3))}$  that
\beno
\varepsilon_7^{\f {q_0} 2}K&\leq& C\sum_{l=1}^K\big(G(u,p,q_0;r,z_l)^{q_0}+H(\pi,p/2,q_0/2;r,z_l)^{q_0/2}\big)\\
&\le& Cr^{\al}\int_{t_0-r^2}^{t_0}\big(\int_{\Om}|u|^pdx\big)^{q_0/p}dt+C(C\theta^{q_0-2})^kr^{\be}\int_{t_0-r^2}^{t_0}\big(\int_{\Om}|\pi|^{\f p2}dx\big)^{q_0/p}dt\\
&&\quad+C\theta^{-2-\frac{3q_0}p+ {q_0} }\sum_{i=0}^{k-1}(C\theta^{q_0-2})^{k-i-1}r^{\al}\int_{t_0-r^2}^{t_0}\big(\int_{\Om}|u|^pdx\big)^{q_0/p}dt\\
&\leq& C\|u\|_{L^{q,\infty}(-1,0;L^p(\R^3))}^{q_0}+C(C\theta^{q_0-2})^k\|\pi\|_{L^{q/2,\infty}(-1,0;L^{p/2}(\Om))}^{q_0/2},
\eeno
where $\Om=\cup_{l=1}^{K}B_{r_0}(x_l), \al=(1-\f 3 p-\f 2{q_0} )q_0, \be=(2-\f 3 p-\f 2{q_0} )$,
and choose $\theta$ such that $C\theta^{2-\frac{4}{q_0}}<\frac12$.
Letting $k\rightarrow\infty$, we infer that
$$K\leq C\varepsilon_7^{-\f {q_0}2}\|u\|_{L^{q,\infty}(-1,0;L^p(\R^3))}^{q_0}.$$

Similarly, for $\frac{q_0}{p}<1$, noting that
$$
\sum_{l=1}^Ka_l^{\f {q_0} p}\le K^{1-\f {q_0} p}\big(\sum_{l=1}^Ka_l\big)^{\f {q_0} p},\quad a_l\ge 0,
$$
we infer that
\beno
\varepsilon_7^{\f {q_0} 2}K
&\leq& C\sum_{l=1}^K\Big(r^{\al}\int_{t_0-r^2}^{t_0}\big(\int_{B_r(x_l)}|u|^pdx\big)^{q_0/p}dt+(C\theta^{q_0-2})^kr^{\be}\int_{t_0-r^2}^{t_0}\big(\int_{B_r(x_l)}|\pi|^{\f p2}dx\big)^{q_0/p}dt\\
&&+C\theta^{q_0-\frac{3q_0}p-2}\sum_{i=0}^{k-1}(C\theta^{q_0-2})^{k-i-1}r^{\al}\int_{t_0-r^2}^{t_0}\big(\int_{B_r(x_l)}|u|^pdx\big)^{q_0/p}dt\Big)\\
&\leq& CK^{1-\f{q_0}p}\Big(r^{\al}\int_{t_0-r^2}^{t_0}(\int_{\Om}|u|^pdx)^{q_0/p}dt+(C\theta^{q_0-2})^kr^{\be}\int_{t_0-r^2}^{t_0}(\int_{\Omega}|\pi|^{p/2}dx)^{q_0/p}dt\Big)\\
&\le& CK^{1-\f{q_0}p}\Big(\|u\|_{L^{q,\infty}(-1,0;L^p(\R^3))}^{q_0}
+(C\theta^{q_0-2})^k\|\pi\|_{L^{q/2,\infty}(-1,0;L^{p/2}(\Om))}^{q_0/2}\Big).
\eeno
Letting $k\rightarrow\infty$, we get
$$K\leq C\varepsilon_7^{-\f p 2}\|u\|_{L^{q,\infty}(-1,0;L^p(\R^3))}^{p}.$$
The proof is completed.\endproof

\section{Appendix}

We first present some estimates of the pressure in terms of some scaling invariant
quantities.

\begin{Lemma}\label{lem:pressure}
Let $(u,\pi)$ be a suitable weak solution of (\ref{eq:NS}) in $Q_1$. Then there hold
\beno
&&H(\pi,2,1;r)\leq C(\frac{\rho}{r})^{\frac32} G(u,4,2;\rho)^2+CH(\pi,1,1;\rho),\\
&&\widetilde{H}(\pi,2,1;r)\leq C(\frac{\rho}{r})^{\f
32}G(u,4,2;\rho)^2+C(\frac{r}{\rho})\widetilde{H}(\pi,1,1;\rho),
\eeno
for any $0<4r<\rho<1$. Here $C$ is a constant independent of $r,\rho$.
\end{Lemma}
\no{\bf Proof.} We write $\pi=\pi_1+\pi_2$ with $\pi_1$ satisfying
\beno
\triangle \pi_1=-\partial_i\partial_j(u_iu_j\zeta),
\eeno
where $\zeta$ is a cut-off function , which equals 1 in $B_{\rho/2}$ and vanishes outside of $B_{\rho}$. Hence,
\beno
\triangle \pi_2=0\quad {\rm in}\quad B_{\rho/2}.
\eeno
By Calder\'{o}n-Zygmund inequality, we have
\beno
\int_{B_{\rho}}|\pi_1|^2dx\leq C\int_{B_{\rho}}|u|^4dx,
\eeno
and using the properties of harmonic function, for $r<\rho/4$
\beno
&&\sup_{x\in B_r}|\pi_2|\leq C\rho^{-3}\int_{B_{\rho/4}}|\pi_2|dx,\\
&&\sup_{x\in B_r}|\pi_2-(\pi_2)_{B_r}|\leq Cr\sup_{x\in B_{\rho/4}}|\nabla \pi_2|
\leq C(\frac{r}{\rho})\rho^{-3}\int_{B_{\rho}}|\pi_2-(\pi_2)_{B_{\rho}}|dx.
\eeno
Then it follows that for $0<r<\rho/4$,
\beno
\int_{B_{r}}|\pi|^2dx
&\leq&\int_{B_{r}}|\pi_1|^2dx+\int_{B_{r}}|\pi_2|^2dx\\
&\leq& C\int_{B_{\rho}}|u|^4dx+Cr^3\rho^{-6}\Big(\int_{B_{\rho}}|\pi|dx\Big)^2,
\eeno
and
\beno
\int_{B_{r}}|\pi-(\pi)_{B_r}|^2dx
&\leq&\int_{B_{r}}|\pi_1-(\pi_1)_{B_r}|^2dx+\int_{B_{r}}|\pi_2-(\pi_2)_{B_r}|^2dx\\
&\leq&C\int_{B_{\rho}}|u|^4dx+Cr^5\rho^{-8}\Big(\int_{B_{\rho}}|\pi-(\pi)_{B_{\rho}}|dx\Big)^2.
\eeno
Integrating with respect to $t$, we get
\beno
\int_{-r^2}^0\Big(\int_{B_{r}}|\pi|^2dx\Big)^{\frac{1}{2}}dt
\leq C\int_{-\rho^2}^0\Big(\int_{B_{\rho}}|u|^4dx\Big)^{\frac{1}{2}}dt
+Cr^{\frac32}\rho^{-3}\int_{-\rho^2}^0\int_{B_{\rho}}|\pi|dxdt,
\eeno
and
\beno
&&\int_{-r^2}^0\Big(\int_{B_{r}}|\pi-(\pi)_{B_r}|^2dx\Big)^{\frac{1}{2}}dt\\
&&\leq C\int_{-\rho^2}^0\Big(\int_{B_{\rho}}|u|^4dx\Big)^{\frac{1}{2}}dt
+r^{\frac52}\rho^{-4}\int_{-\rho^2}^0\int_{B_{\rho}}|\pi-(\pi)_{B_{\rho}}|dxdt.
\eeno
The proof is completed.\endproof

The same proof also yields that for any $0<4r<\rho<1$,
\ben\label{eq:pressure1}
H(\pi,p/2,q/2;r)\leq C\big(\frac{\rho}{r}\big)^{\frac{4}{q}+\frac6p-2}G(u,p,q;\rho)^2
+C\big(\frac{r}{\rho}\big)^{2-\frac{4}{q}}H(\pi,1,q/2;\rho),
\een
where $p>2, q\ge 2$.
Similarly, one can show that (see also \cite{Seregin-JMS})
\ben\label{eq:pressure}
\widetilde{D}(\pi,r)\leq C\big((\frac{r}{\rho})^{5/2}\widetilde{D}(\pi,\rho)+(\frac{\rho}{r})^{2}C(u,\rho)\big),
\een
for any $0<4r<\rho<1.$

The following lemma gives a bound of local scaling invariant energy, see also \cite{GKT} and \cite{ZS}.

\begin{Lemma}\label{lem:invariant}
Let $(u,\pi)$ be a suitable weak solution of (\ref{eq:NS}) in $Q_1$. If
\beno
&&G(u,p,q;r)\leq M\quad\textrm{ with }\quad1\le\frac3p+\frac2q<2, 1<q\leq \infty\quad \textrm{or}\\
&&H(\nabla u, p, q;r)\leq M\quad\textrm{ with }\quad 2\le \frac3{p}+\frac2{q}<3, 1<p\le\infty,
\eeno
for any $0<r<1$, then there holds for $0<r<1/2$
\beno
A(u,r)+E(u,r)+D(\pi,r)\leq C(p,q,M)\big(r^{1/2}\big(C(u,1)+D(\pi,1)\big)+1\big).
\eeno
\end{Lemma}

\no{\bf Proof.}\,First of all, we assume $G(u,p,q;r)\leq M$ and moreover,
\beno
\f32<\frac3p+\frac2q<2,\quad \frac3p+\frac3q\geq 2,\quad \frac4p+\frac2q\geq 2,\quad p,q<\infty.
\eeno
Otherwise, we can choose $(p_1,q_1)$ satisfying the above condition and by H\"{o}lder inequality,
\beno
G(u,p_1,q_1;r)\le CG(u,p,q;r).
\eeno

By H\"{o}lder inequality and Sobolev inequality, we get
\beno
\int_{B_r}|u|^3dx &= &\int_{B_r}|u|^{3\alpha+3\beta+3-3\alpha-3\beta}dx\\
&\leq & \big(\int_{B_r}|u|^2dx\big)^{3\alpha/2}\big(\int_{B_r}|u|^6dx\big)^{\beta/2}
\big(\int_{B_r}|u|^pdx\big)^{{(3-3\alpha-3\beta)}/p}\\
&\leq & C\big(\int_{B_r}|u|^2dx\big)^{3\alpha/2}\big(\int_{B_r}|\nabla u|^2+|u|^2dx\big)^{3\beta/2}
\big(\int_{B_r}|u|^pdx\big)^{{(3-3\alpha-3\beta)}/p},
\eeno
where  $\alpha, \beta\geq 0$ are chosen such that
\beno
\frac13=\frac{\alpha}{2}+\frac{\beta}{6}+\frac{1-\alpha-\beta}{p},\quad
1=\frac{3\beta}{2}+\frac{3-3\alpha-3\beta}{q}.
\eeno
That is,
\beno
\alpha=\frac{2(\frac3p+\frac3q-2)}{3(\frac6p+\frac4q-3)},\quad \beta=\frac{\frac4p+\frac2q-2}{\frac6p+\frac4q-3}.
\eeno
Integrating with respect to time, we get
\beno
\int_{Q_r}|u|^3dxdt
&\leq& C\big(\sup_{-r^2<t<0}\int_{B_r}|u|^2dx\big)^{\frac{3\alpha}{2}}
\big(\int_{Q_r}|\nabla u|^2+|u|^2dxdt\big)^{\f {3\beta}{2}}\\
&&\quad\times \Big(\int_{-r^2}^0\big(\int_{B_r}|u|^pdx\big)^{\f q p}dt\Big)^{\f {3-3\alpha-3\beta} q},
\eeno
this means that
\beno
C(u,r)\leq C\big(A(u,r)+E(u,r)\big)^{\f {3\alpha+3\beta}{2}}G(u,p,q;r)^{3-3\alpha-3\beta}.
\eeno
Set $\frac3p+\frac2q=2-\delta$ with $0\leq\delta<1/2$.
Then $\f {3\alpha+3\beta}{2}=\frac32-\frac{1}{2(\frac6p+\frac4q-3)}=\frac{1-3\delta}{1-2\delta}$ and
\ben\label{eq:7.9}
C(u,r)\leq C\big(A(u,r)+E(u,r)\big)^{\frac{1-3\delta}{1-2\delta}}G(u,p,q;r)^{\frac{1}{1-2\delta}}.
\een
By the assumption, we get
\beno
C(u,r)\leq C(p,q,M)\big(A(u,r)+E(u,r)\big)^{\frac{1-3\delta}{1-2\delta}}.
\eeno
Using the local energy inequality and (\ref{eq:pressure1}), we deduce that
\beno
&&A(u,r)+E(u,r)\leq C\big(C(u,2r)^{2/3}+C(u,2r)+C(u,2r)^{1/3}D(\pi,2r)^{2/3}\big),\\
&&D(\pi,r)\leq C\big((\frac{r}{\rho})D(\pi,\rho)+(\frac{\rho}{r})^{2}C(u,\rho)\big)\quad \textrm{for}\quad 0<4r<\rho<1.
\eeno
Set $F(r)=A(u,r)+E(u,r)+D(\pi,r)$. It follows from the above three inequalities that
\beno
F(r)&\leq& C\big(1+C(u,2r)+D(\pi,2r)\big)\\
&\leq& C+C(\frac{r}{\rho})F(\rho)+
C(p,q,M)\big((\frac{\rho}{r})^{2}+(\frac{\rho}{r})^{\frac{1-3\delta}{1-2\delta}}\big)\big(A(u,\rho)+E(u,\rho)\big)^{\frac{1-3\delta}{1-2\delta}}\\
&\leq& C+C(\frac{r}{\rho})F(\rho)+C(p,q,M,\frac{\rho}{r})
\eeno
for $0<8r<\rho<1$. By the standard iteration and local energy inequality, we deduce that
\beno
F(r)&\leq& C(p,q,M)\big(r^{1/2}(A(u,1/2)+E(u,1/2)+D(\pi,1))+1\big)\\
&\leq& C(p,q,M)\big(r^{1/2}(C(u,1)+D(\pi,1))+1\big).
\eeno

Now Let us assume that $H(\nabla u, p, q;r)\leq M$ and $\f52<\frac3{p}+\frac2{q}<3, p<3$.
General case can be reduced to this case as above. Similarly, we have
\beno
&&\int_{Q_r}|u-u_{B_r}|^3dxdt\\
&&\leq C\big(\sup_{-r^2<t<0}\int_{B_r}|u|^2dx\big)^{\frac{3\alpha}{2}}
\big(\int_{Q_r}|\nabla u|^2dxdt\big)^{\f {3\beta}{2}}\\
&&\quad\times\Big(\int_{-r^2}^0\big(\int_{B_r}|u-u_{B_r}|^{\frac{3p}{3-p}}dx\big)^{\f {q(3-p)}{3p} }dt\Big)^{\f {3-3\alpha-3\beta} {q}}\\
&&\leq C\big(\sup_{-r^2<t<0}\int_{B_r}|u|^2dx\big)^{\frac{3\alpha}{2}}
\big(\int_{Q_r}|\nabla u|^2dxdt\big)^{\f {3\beta}{2}}
\Big(\int_{-r^2}^0\big(\int_{B_r}|\nabla u|^{p}dx\big)^{\f {q}{p} }dt\Big)^{\f {3-3\alpha-3\beta} {q}},
\eeno
where $\alpha+\beta=1-\frac{1}{3(\frac{6}{p}+\frac4{q}-5)}$. Let $\frac{3}{p}+\frac2{q}=3-\delta_0$ with
$0\leq \delta_0<\frac12$, then
\beno
\widetilde{C}(u,r)&\leq& C\big(A(u,r)+E(u,r)\big)^{\frac{1-3\delta_0}{1-2\delta_0}}H(\nabla u,p,q;r)^{\frac{1}{1-2\delta_0}}\\
&\leq& C(p,q,M)\big(A(u,r)+E(u,r)\big)^{\frac{1-3\delta_0}{1-2\delta_0}}.
\eeno
Note that
\beno
C(u,r)\leq C\big((\frac{r}{\rho})C(u,\rho)+(\frac{\rho}{r})^{2}\tilde{C}(u,\rho)\big),
\eeno
and
\beno
&&A(u,r)+E(u,r)\leq C\big(C(u,2r)^{2/3}+C(u,2r)+C(u,2r)^{1/3}D(\pi,2r)^{2/3}\big),\\
&&D(\pi,r)\leq C\big((\frac{r}{\rho})D(\pi,\rho)+(\frac{\rho}{r})^{2}\widetilde{C}(u,\rho)\big),
\eeno
for $0<4r<\rho<1$. Let $F(r)=A(u,r)+E(u,r)+C(u,r)+D(\pi,r)$. Then we have
\beno
F(r)&\leq& C\big(1+C(u,2r)+D(\pi,2r)\big)\\
&\leq& C\big(1+(\frac{r}{\rho})C(u,\rho)+(\frac{\rho}{r})^{2}\widetilde{C}(u,\rho)
+(\frac{r}{\rho})D(\pi,\rho)\big)\\
&\leq& C\big(1+(\frac{r}{\rho})F(\rho)+(\frac{\rho}{r})^{2}F(\rho)^{\frac{1-3\delta_0}{1-2\delta_0}}\big)\\
&\leq& C+C(\frac{r}{\rho})F(\rho)+C(p,q,M,\frac{\rho}{r}),
\eeno
which implies the required result.\endproof

\bigskip

\noindent {\bf Acknowledgments.}
Both authors thank helpful discussions with Professors Gang Tian and Liqun Zhang.
Zhifei Zhang is partly supported by NSF of China under Grant 10990013 and 11071007.

\end{document}